\documentclass[12pt]{amsart}
\usepackage{latexsym,amssymb,verbatim,multirow,graphicx,epsfig,enumerate, paralist,hyperref}
\usepackage{xcolor}
\usepackage[normalem]{ulem}

\theoremstyle{plain}
   \newtheorem{theorem}{Theorem}
   
   \newtheorem{prop}[theorem]{Proposition}

   \newtheorem*{theorem*}{Theorem}
   
\theoremstyle{definition}
   \newtheorem{definition}[theorem]{Definition}

\numberwithin{equation}{section}

\newcommand\FF{{\mathbb{F}}}

\newcommand\GL{{\mathrm{GL}}}

\newcommand{\GLnFp}{\GL_n(\Fp)}

\newcommand{\Fp}{\FF_p}

\title{On the Brunnian Conjecture}
\author{Jean-Yves Degos}

\begin{document}

\begin{abstract}
We prove the Conjecture $B(n,p,P)$ of \cite[Conjecture 3.6 p.~62]{Degos}, which is renamed here as ``the Brunnian Conjecture'', in the case where $p \geq 5$ and $n$ is not a multiple of $p-1$.
\end{abstract}

\maketitle

The goal of this paper is to give a partial proof of the Conjecture $B(n,p,P)$ (\cite[Conjecture 3.6 p.~62]{Degos}). We shall refer to it as the ``Brunnian Conjecture'' from now on. Our strategy consists in using Nick Gill's proof of the Conjecture $A(n,p,P)$ (\cite[Conjecture 3.4 p.~62]{Degos}) given in \cite[Theorem 2 p.~230]{Gill}, corrected by the paper of Joel Brewster Lewis, see (\cite[Theorem 7 p.~6 and Correction 3.3 p.~10]{Lewis}).

Our main result is:

\begin{theorem*} Let $p$ a prime number, and $n \geq 3$ an integer, and
\[ f(X) = X^n + a_{n-1} X^{n-1} + \cdots + a_1 X + a_0 \in \Fp[X] \]
a primitive polynomial of degree $n$. Let $C=C_f$ denote the companion matrix of $f$, namely:
\[ C = 
\left[
\begin{array}{ccccc}
0 &  0 & \ldots   & 0 & -a_0 \\
1 &  0 & \ldots   & 0 & -a_1 \\
0 &  \ddots & \ddots   & \vdots & \vdots \\
\vdots &  \ddots & \ddots   & 0 & -a_{n-2} \\
0 &  \ldots & 0  & 1 & -a_{n-1}
\end{array}
\right].
\]

Let  $G$ denote the companion matrix of the polynomial $X^n-1$, and let:
\[ \begin{array}{rcl}
G_1 & := & C \,, \\
G_{k+1} & := & G G_k G^{-1}  \text{for all } 1 \leq k \leq n-1.
\end{array} 
\]
Then, if $p \geq 5$ and $n$ is not divisible by $p-1$, then 
\[ \langle G_1, G_2, \ldots, G_n \rangle = \GLnFp \,. \] 
\end{theorem*}

\section{Background and strategy}

This story started with René Guitart's works on Klein's group (defined as the autormorphisms group of the Klein's quartic with 168 elements), described as a borromean objet (\cite[section 3 p. 148]{Guitart}), which is nothing else that $\text{GL}_3(\mathbb{F}_2)$. In \cite{Degos}, we extended René Guitart's notion of borromean group by introducing the notions of $n$-cyclable groups (\cite[Definition 1.1 p. 58]{Degos}) and $n$-brunnian groups (\cite[Definition 1.2 and Definition 1.3 p. 58]{Degos}. The so-called ``Brunnian conjecture'' refers to the following statement: if $p$ is a prime number, and $n \geq 2$, then the group $\GLnFp$ can be generated by the companion matrix of a primitive polynomial $P(X) \in \Fp[X]$ of degree $n$, and the conjugates by $G^k$ for $1 \leq k \leq n-1$ of the companion matrix of polynomial $X^n-1$. This conjecture is correlated to this other one:  the group $\GLnFp$ can be generated by the companion matrix of a primitive polynomial $P(X) \in \Fp[X]$ of degree $n$, denoted by $C$, and the companion matrix of the polynomial $X^n-1$, denoted by $G$. Thus:
\[
\GLnFp = \langle C, G \rangle .
\]

Then, Nick Gill, using the classification of finite simple groups, gave a proof of the latter statement, as a consequence of the following slightly more general result (\cite[Theorem 2 p. 230]{Gill}): if $\FF$ is a finite field of order $q$ and $f, g$ are distinct polynomials of degree $n$ such that $f$ is primitive, and if the constant term of $g$ is non-zero, then:
\[
\langle C_f, C_g \rangle = \GLnFp.
\]

Consequently, we can infer from Nick Gill's theroem a quite simple way to prove the Brunnian conjecture: it sufficient to build a companion matrix $K$ (of a polynomial $g$ of degree $n$) with the following properties:
\begin{itemize}
    \item $K \in \langle G_1, G_2, \ldots, G_n \rangle$;
    \item $g(0) \neq 0$;
    \item $K \neq G_1$.
\end{itemize}
Thus, as $\langle G_1, K \rangle \subset \langle G_1, G_2, \ldots, G_n \rangle \subset \GLnFp$, whenever  $\langle G_1, K \rangle = \GLnFp$, we are done!

In this paper, our main contribution is to build such a matrix $K$ in a generic way... Our machinery works whenever $-a0 \neq 1$ and $(-a_0)^n \neq 1$, which is guaranteed by the assumptions $p \geq 5$ and $n$ not divisible par $p-1$. Morever, we must assume that $n \geq 3$, as Joel Brewster pointed it out.

\section{Statements and proofs}


\begin{definition} With the notations of the main theorem, assuming $n \geq 3$, we define the matrix $K$ by:
\[
K := G_n G_{n_1} \cdots G_2 G_1 \,.
\]
\end{definition}

\begin{prop}
We have $K = C (G^{-1}C)^n$.
\end{prop}

\begin{proof}
We start by showing, thanks to a recurrence-based reasoning, the following result:
\[ G_k G_{k-1} \ldots G_1 = G^{k-1} C (G^{-1} C)^{k-1} \text{ for all } 2 \leq k \leq n \,.\] 
Indeed, if $k=2$, this equality is equivalent to: 
\[ G_2 G_1 = G^{2-1} C (G^{-1} C)^{2-1}~; \]
therefore, it is true because the right hand side equals $G C G^{-1} C = G_2 G_1$. We assume the recurrence hypothesis at rank $k$, and we show if at rank $k+1$. We have:
\begin{eqnarray*}
    G_{k+1} G_k G_{k-1} \cdots G_1 & = & G_{k+1} \cdot G^{k-1} C (G^{-1} C)^{k-1} \\
    & = & G^k C G^{-k} \cdot G^{k-1} C (G^{-1} C)^{k-1} \\
    & = & G^k \cdot C \cdot G^{-k} \cdot G^{k-1} C \cdot (G^{-1} C)^{k-1} \\
    & = & G^k \cdot C \cdot G^{-1} C \cdot (G^{-1} C)^{k-1} \\
    & = & G^k \cdot C \cdot (G^{-1} C)^{k} \,.\\
\end{eqnarray*}
We then get: $ G_n G_{n-1} \ldots G_1 = G^{n-1} C \bigl( G^{-1} \bigr)^{n-1} C$.  We conclude that:
\begin{eqnarray*}
G_1 G_n G_{n-1} \ldots G_1 & = & G_1 G^{n-1} C ( G^{-1})^{n-1} \,,  \\
 & = & C G^{-1} C ( G^{-1})^{n-1} \,, \\
 & = & C (G^{-1} C)^n \,,
\end{eqnarray*}
QED.
\end{proof}


\begin{prop} 
\label{propgminus1c}
Let $E_k$ denote the $k^\text{th}$ vector of the canonical basis, for $1 \leq k \leq n$. Then we get:
\begin{enumerate}
\item $G^{-1}C.E_k = E_k$ for $1  \leq k \leq n-1;$
\item $G^{-1}C.E_n = -a_0 E_n + \sum_{j=0}^{n-1} -a_j E_j,$
\end{enumerate}
id est~:
$$ G^{-1}C = 
\left[
\begin{array}{ccccc}
1 &  0 & \ldots   & 0 & -a_1 \\
0 &  1 & \ldots   & 0 & -a_2 \\
\vdots &  \vdots & \ddots   & \vdots & \vdots \\
0 &  0 & \ldots   & 1 & -a_{n-1} \\
0 &  0 & \ldots   & 0 & -a_0
\end{array}
\right] \,.
$$
\end{prop}

\begin{proof}
The first statement is a consequence of the fact that the product of the matrix $G$ (respectively $C$) by the column vector $E_k$ is equal to the column vector $E_{k+1}$ for $1 \leq k \leq n-1$. Thus, $G^{-1}C.E_k = G^{-1}.CE_k = G^{-1}. E_{k+1} = E_k$.

To show the second statement, we proceed as follows:
\begin{eqnarray*}
    G^{-1}C.E_n & = & G^{-1}.\left( \sum_{i=1}^n -a_{i-1} E_i \right) \\
    & = & \sum_{i=1}^n -a_{i-1} G^{-1}.E_i \\
    & = & \sum_{i=1}^n -a_{i-1} \left\{ \begin{array}{rcl} E_{i-1} & \text{ si } 2 \leq i \leq n \\ E_n & \text{ si } & i = 1 \end{array} \right. \\
    & = & -a_0 G^{-1}E_1 + \sum_{i = 2}^n -a_{i-1} G^{-1} E_i \\
    & = & -a_0 E_n + \sum_{i = 2}^n -a_{i-1} E_{i-1} \\
    & = & -a_0 E_n + \sum_{j = 1}^n -a_{j} E_{j} \\
\end{eqnarray*}
thanks to the change of variables $j = i-1$.
\end{proof}


\begin{prop}
We assume that $-a_0 \neq 1$. Then be the matrix $G^{-1} C$ can be written as $Q^{-1} G^{-1}C Q = D$, where $Q \in \GLnFp$ and $D$, diagonal matrices, express like~:
$$ Q = 
\left[
\begin{array}{ccccc}
1 &  0 & \ldots   & 0 & \frac{a_1}{1+a_0} \\
0 &  1 & \ldots   & 0 & \frac{a_2}{1+a_0} \\
\vdots &  \vdots & \ddots   & \vdots & \vdots \\
0 &  0 & \ldots   & 1 & \frac{a_{n-1}}{1+a_0} \\
0 &  0 & \ldots   & 0 & 1
\end{array}
\right] \text{ and } D = 
\left[
\begin{array}{ccccc}
1 &  0 & \ldots   & 0 & 0 \\
0 &  1 & \ldots   & 0 & 0 \\
\vdots &  \vdots & \ddots   & \vdots & \vdots \\
0 &  0 & \ldots   & 1 & 0 \\
0 &  0 & \ldots   & 0 & -a_0
\end{array}
\right] \,,
$$
We deduce that: $K = G^{n-1} (G^{-1}C)^{n} = C (Q D Q^{-1})^n = C Q D^n Q^{-1} \,.$
\end{prop}

\begin{proof}

The chacteristic polynomial of the matrix $G^{-1} C$ is $\chi_{G^{-1} C}(X) = (X-1)^{n-1} (X+a_0)$. Therefore, the eigenvalues are $1$, of order $n-1$ and $-a_0$, or order $1$. From Propostion~\ref{propgminus1c}, (1), page~\pageref{propgminus1c}, we can infer that  $Q$ will satisfy $Q.E_k = E_k$ for $1 \leq k \leq n-1$. Consequently, we must find a column vector $V$ such that $G^{-1} C V = -a_0 V$, as eigenvector associated to the eigenvalue $-a_0$.

Thefore, we write~:
\[
V := \left( \begin{array}{c}
    v_1  \\
    v_2 \\
    \vdots \\
    v_n
\end{array}\right), \text{ hence } 
G^{-1} C V = \left( \begin{array}{c}
    v_1 - a_1 v_n \\
    v_2 - a_2 v_n \\
    \vdots \\
    v_{n-1} - a_{n-1} v_n \\
    -a_0 v_n 
\end{array}\right), \text{ and } 
\]
\[
-a_0 V = \left( \begin{array}{c}
    - a_0 v_1 \\
    - a_0 v_2 \\
    \vdots \\
    -a_0 v_n
\end{array}\right)
\]
Solving this Cramer linear system in $v_1, v_2, \ldots, v_n$, we get:
\begin{enumerate}
\item $v_n \in \mathbb{F}_p \,;$
\item $v_i = \frac{a_i v_n }{1+a_0}$ for $1 \leq i \leq n-1 \,.$
\end{enumerate}

We now just have to choose $v_n = 1$ to get the previously announced matrix $Q$.
\end{proof}


\begin{prop}
If $- a_0 \neq 1$ and $(-a_0)^n \neq 1$, the matrix 
$$K:= G_n G_{n-1} \cdots G_2 G_1$$
is a companion matrix, different from $C$.
\end{prop}

\begin{proof}

Under the action of $Q$ and $D$, the columun vectors $E_k$ for $1 \leq k \leq n-1$ remain unchanged. We thus have:
\begin{eqnarray*}
    K.E_k & = & C Q D^n Q^{-1}.E_k \\
    & = & C.E_k \\
    & = & E_{k+1} \,.
\end{eqnarray*}
Therefore, we must evaluate $K.E_n = C Q D^n Q^{-1}.E_n$

Let $U := Q^{-1}E_n$. Then $E_n = Q U$. To find the valeur of $U$, it is sufficient to solve the equation. However, if
\[
Q = 
\left[
\begin{array}{ccccc}
1 &  0 & \ldots   & 0 & q_1 \\
0 &  1 & \ldots   & 0 & q_2 \\
\vdots &  \vdots & \ddots   & \vdots & \vdots \\
0 &  0 & \ldots   & 1 & q_{n-1} \\
0 &  0 & \ldots   & 0 & q_n
\end{array}
\right] \text{ and } U = 
\left[
\begin{array}{c}
u_1 \\
u_2 \\
\vdots \\
u_{n-1} \\
u_n
\end{array}
\right] \,, 
\]
then
\[
QU = 
\left[
\begin{array}{c}
u_1 + q_1 u_n \\
u_2 + q_2 u_n\\
\vdots \\
u_{n-1} + q_{n-1} u_n \\
q_n u_n
\end{array}
\right] \text{ equalized to } E_n = 
\left[
\begin{array}{c}
0 \\
0 \\
\vdots \\
0 \\
1
\end{array}
\right] 
\]
enables to get~: $u_n = \frac{1}{q_n}  = 1$ (because $q_n$ is a notation for the coefficient in row $n$ and columun $n$ of $Q$, which equals $1$), and:
\[  u_k = - q_k u_n = - q_k = \frac{-a_k}{1+a_0} \,.\]

Thus, we have:
\[
U = 
\left[
\begin{array}{c}
\frac{-a_1}{1+a_0} \\
\frac{-a_2}{1+a_0}  \\
\vdots \\
\frac{-a_{n-1}}{1+a_0} \\
1
\end{array}
\right] \,.
\]

Now, let $W:= D^n U$. Then:
\[
W = 
\left[
\begin{array}{c}
\frac{-a_1}{1+a_0} \\
\frac{-a_2}{1+a_0}  \\
\vdots \\
\frac{-a_{n-1}}{1+a_0} \\
(-a_0)^n
\end{array}
\right] \,,
\]
and $CQD^n U = CQ W$. Now, the following holds~:
\[
QW = 
\left[
\begin{array}{c}
0 \\
0  \\
\vdots \\
0 \\
(-a_0)^n
\end{array}
\right] \,,
\]
To summarize, we get:
\begin{eqnarray*}
    K.E_n & = & CQ D^n Q^{-1}.E_n \\
    & = & CQ D^n U \\
    & = & CQ WW \\
    & = & C (-a_0)^n E_n \\
    & = & (-a_0)^n C E_n \\
    & = & \sum_{i = 1}^n -a_{i-1} (-a_0)^n E_i \,.
\end{eqnarray*}

Therefore, the matrix $K$ is a companion matrix, different from $C$, because we assumed that $(-a_0)^n \neq 1$.    
\end{proof}


\begin{theorem*}[Main result]
Let $p \geq 5$ a prime number. Then $-a_0 \neq 1$. Moreover, if $n$ is not divisible by $p-1$ then $(-a_0)^n \neq 1$. In this case, the Brunnian conjecture $B(n,p,P)$ is true.
\end{theorem*}

\begin{proof}

As the polynomial $f$ is primitive, $(-1)^n a_0$ generates the multiplicative group $\mathbb{F}_p^\times$, so if $p \geq 5$, the equality $-a_0 = 1$ can not occur. Indeed, if $-a_0 = 1$ held, we would have $a_0 = -1$, then $(-1)^n a_0 = (-1)^{n+1}$ would be equal to $-1$ or $1$, which is never a generator of $\mathbb{F}_p^\times$ when $p \geq 5$.

If $n$ if not divisible by $p-1$, let us assume that $(-a_0)^n = 1$ in $\Fp^\times$. We can write $n = (p-1)q + r$, where $1 \leq r < p-1$ is an integer. Then $(-a_0)^r = 1$, because according to our assumption, $1 = (-a_0)^n = (-a_0)^{(p-1)q} (-a_0)^r$ and $(-a_0)^{p-1}=1$ (because $-a_0 \in \Fp^\times$). Now, it turns out that $(-1)^n a_0 = (-1)^{r-1}(-a_0)$ has an order which divides $r < p-1$. This is a contradiction with the fact that $(-1)^n a_0$ is a generator of $\Fp^\times$, the order of which is $p-1$. Consequently, if $n$ if not divisible by $p-1$, then $(-a_0)^n \neq 1$.
\end{proof}

\section{Further remarks}

In \cite{Degos}, we gave some ``experimental checkings'', and we were unable to conclude whether the conjecture was true for certain couples $(n,p)$, due to a lack of power of computations and/or high complexity of calculus. Thanks to our new method, some answers can be given. We sum up them in the following tables. If ``(False)'' is written in last column, it means that the conjecture can be true, but it is not a consequence of the main result of this paper.

\begin{center}

\pagebreak

Primitive polynomials over $\FF_5$

\scriptsize

\begin{tabular}{|c|c|c|c|c|} \hline
$n$ & $p$ & $P(x)$ & $B(n,p,P)$ \\ \hline
$7$ & $5$ & $x^7+x^6+2$ & True \\
$8$ & $5$ & $x^8+x^5+x^3+3$ & (False) \\
$9$ & $5$ & $x^9+x^7+x^6+3$ & True \\
$10$ & $5$ & $x^{10}+x^9+x^7+3$ & True \\
$11$ & $5$ & $x^{11}+x^{10}+2$ & True \\
$12$ & $5$ & $x^{12}+x^7+x^4+3$ & (False) \\
$13$ & $5$ & $x^{13}+4x^2+3x+3$ & True \\
$14$ & $5$ & $x^{14}+x^7+4x^5+4x^4+2x^3+3x^2+x+2$ & True \\
$15$ & $5$ & $x^{15}+2x^5+4x^5+3x^3+3x^2+4x+3$ & True \\
$16$ & $5$ & $x^{16}+x^7+4x^6+4x^5+4x^4+4x^3+4x^2+x+2$ & (False) \\
$17$ & $5$ & $x^{17}+3x^2+2x+3$ & True \\
$18$ & $5$ & $x^{18}+x^{12}+x^{11}+x^{10}+x^9+2x^8+2x^6+x^5+2x^3+2x^2+2$ & True \\
\hline
\end{tabular}

\normalsize

\vspace{1em}

Primitive polynomials over $\FF_7$

\begin{tabular}{|c|c|c|c|c|} \hline
$n$ & $p$ & $P(x)$ & $B(n,p,P)$ \\ \hline
$6$ & $7$ & $x^6+x^5+x^4+3$ & (False) \\
$7$ & $7$ & $x^7+x^5+4$ & True \\
$8$ & $7$ & $x^8+x^7+3$ & True \\
$9$ & $7$ & $x^9+x^8+x^3+2$ & True \\
$10$ & $7$ & $x^{10}+x^9+x^8+3$ & True \\
\hline
\end{tabular}

\vspace{1em}

Primitive polynomials over $\FF_{11}$

\begin{tabular}{|c|c|c|c|c|} \hline
$n$ & $p$ & $P(x)$ & $B(n,p,P)$ \\ \hline
$5$ & $11$ & $x^5+x^4+x^3+3$ & True \\
$6$ & $11$ & $x^6+x^5+x+7$ & True \\
$7$ & $11$ & $x^7+x^6+4$ & True \\
$8$ & $11$ & $x^8+x^7+x^6+7$ & True\\
\hline
\end{tabular}

\vspace{1em}

Primitive polynomials over $\FF_{13}$

\begin{tabular}{|c|c|c|c|c|} \hline
$n$ & $p$ & $P(x)$ & $B(n,p,P)$ \\ \hline
$5$ & $13$ & $x^5+x^4+x^3+6$ & True \\
$6$ & $13$ & $x^6+x^5+x^3+6$ & True \\
$7$ & $13$ & $x^7+x^4+2$ & True \\
$8$ & $13$ & $x^8+x^7+x^6+11$ & True \\
\hline
\end{tabular}

\vspace{1em}

Primitive polynomials over $\FF_{17}$

\begin{tabular}{|c|c|c|c|c|} \hline
$n$ & $p$ & $P(x)$ & $B(n,p,P)$ \\ \hline
$4$ & $17$ & $x^4+x^3+5$ & True \\
$5$ & $17$ & $x^5+x^4+5$ & True \\
$6$ & $17$ & $x^6+x^5+3$ & True \\
$7$ & $17$ & $x^7+x^6+7$ & True \\
\hline
\end{tabular}

\pagebreak

Primitive polynomials over $\FF_{19}$

\begin{tabular}{|c|c|c|c|c|} \hline
$n$ & $p$ & $P(x)$ & $B(n,p,P)$ \\ \hline
$4$ & $19$ & $x^4+x^3+2$ & True \\
$5$ & $19$ & $x^5+x^4+5$ & True \\
$6$ & $19$ & $x^6+x^5+15$ & True \\
$7$ & $19$ & $x^7+x^6+5$ & True \\
\hline
\end{tabular}

\vspace{1em}

Primitive polynomials over $\FF_{23}$

\begin{tabular}{|c|c|c|c|c|} \hline
$n$ & $p$ & $P(x)$ & $B(n,p,P)$ \\ \hline
$4$ & $23$ & $x^4+x^3+20$ & True  \\
$5$ & $23$ & $x^5+x^4+6$ & True \\
$6$ & $23$ & $x^6+x^5+7$ & True \\
\hline
\end{tabular}

\vspace{1em}

Primitive polynomials over $\FF_{29}$

\begin{tabular}{|c|c|c|c|c|} \hline
$n$ & $p$ & $P(x)$ & $B(n,p,P)$ \\ \hline
$4$ & $29$ & $x^4+x^3+2$ & True \\
$5$ & $29$ & $x^5+x^4+2$ & True \\
$6$ & $29$ & $x^6+x^5+11$ & True \\
\hline
\end{tabular}

Primitive polynomials over $\FF_{31}$

\begin{tabular}{|c|c|c|c|c|} \hline
$n$ & $p$ & $P(x)$ & $B(n,p,P)$ \\ \hline
$4$ & $31$ & $x^4+x^3+13$ & True \\
$5$ & $31$ & $x^5+x^4+10$ & True \\
$6$ & $31$ & $x^6+x^5+12$ & True \\
\hline
\end{tabular}

\vspace{1em}

Primitive polynomials over $\FF_{37}$

\begin{tabular}{|c|c|c|c|c|} \hline
$n$ & $p$ & $P(x)$ & $B(n,p,P)$ \\ \hline
$3$ & $37$ & $x^3+x^2+17$ & True \\
$4$ & $37$ & $x^4+x^3+22$ & True \\
$5$ & $37$ & $x^5+x^4+2$ & True \\
\hline
\end{tabular}

\vspace{1em}

Primitive polynomials over $\FF_{41}$

\begin{tabular}{|c|c|c|c|c|} \hline
$n$ & $p$ & $P(x)$ & $B(n,p,P)$ \\ \hline
$3$ & $41$ & $x^3+x^2+11$ & True \\
$4$ & $41$ & $x^4+x^3+26$ & True \\
$5$ & $41$ & $x^5+x^4+11$ & True \\
\hline
\end{tabular}

\vspace{1em}

Primitive polynomials over $\FF_{43}$

\begin{tabular}{|c|c|c|c|c|} \hline
$n$ & $p$ & $P(x)$ & $B(n,p,P)$ \\ \hline
$3$ & $43$ & $x^3+x^2+9$ & True \\
$4$ & $43$ & $x^4+x+20$ & True \\
$5$ & $43$ & $x^5+x^4+9$ & True \\
\hline
\end{tabular}

\pagebreak

Primitive polynomials over $\FF_{47}$

\begin{tabular}{|c|c|c|c|c|} \hline
$n$ & $p$ & $P(x)$ & $B(n,p,P)$ \\ \hline
$3$ & $47$ & $x^3+x^2+2$ & True \\
$4$ & $47$ & $x^4+x^3+5$ & True \\
$5$ & $47$ & $x^5+x^4+6$ & True \\
\hline
\end{tabular}

\vspace{1em}

Primitive polynomials over $\FF_{53}$

\begin{tabular}{|c|c|c|c|c|} \hline
$n$ & $p$ & $P(x)$ & $B(n,p,P)$ \\ \hline
$3$ & $53$ & $x^3+x^2+2$ & True \\
$4$ & $53$ & $x^4+x^3+2$ & True \\
$5$ & $53$ & $x^5+x^4+12$ & True \\
\hline
\end{tabular}

\vspace{1em}

Primitive polynomials over $\FF_{59}$

\begin{tabular}{|c|c|c|c|c|} \hline
$n$ & $p$ & $P(x)$ & $B(n,p,P)$ \\ \hline
$3$ & $59$ & $x^3+x^2+9$ & True \\
$4$ & $59$ & $x^4+x^3+18$ & True \\
$5$ & $59$ & $x^5+x^4+4$ & True \\
\hline
\end{tabular}

\vspace{1em}

Primitive polynomials over $\FF_{61}$

\begin{tabular}{|c|c|c|c|c|} \hline
$n$ & $p$ & $P(x)$ & $B(n,p,P)$ \\ \hline
$3$ & $61$ & $x^3+x^2+6$ & True \\
$4$ & $61$ & $x^4+x^3+17$ & True \\
$5$ & $61$ & $x^5+x^4+55$ & True \\
\hline
\end{tabular}

\vspace{1em}

Primitive polynomials over $\FF_{67}$

\begin{tabular}{|c|c|c|c|c|} \hline
$n$ & $p$ & $P(x)$ & $B(n,p,P)$ \\ \hline
$3$ & $67$ & $x^3+x^2+6$ & True \\
$4$ & $67$ & $x^4+x^3+12$ & True \\
\hline
\end{tabular}

\vspace{1em}

Primitive polynomials over $\FF_{71}$

\begin{tabular}{|c|c|c|c|c|} \hline
$n$ & $p$ & $P(x)$ & $B(n,p,P)$ \\ \hline
$3$ & $71$ & $x^3+x^2+8$ & True  \\
$4$ & $71$ & $x^4+x^3+13$ & True \\
\hline
\end{tabular}

\vspace{1em}

Primitive polynomials over $\FF_{73}$

\begin{tabular}{|c|c|c|c|c|} \hline
$n$ & $p$ & $P(x)$ & $B(n,p,P)$ \\ \hline
$3$ & $73$ & $x^3+x^2+5$ & True  \\
$4$ & $73$ & $x^4+x^3+33$ & True  \\
\hline
\end{tabular}

\vspace{1em}

Primitive polynomials over $\FF_{79}$

\begin{tabular}{|c|c|c|c|c|} \hline
$n$ & $p$ & $P(x)$ & $B(n,p,P)$ \\ \hline
$3$ & $79$ & $x^3+x^2+2$ & True \\
$4$ & $79$ & $x^4+x^3+7$ & True \\
\hline
\end{tabular}

\vspace{1em}

Primitive polynomials over $\FF_{83}$

\begin{tabular}{|c|c|c|c|c|} \hline
$n$ & $p$ & $P(x)$ & $B(n,p,P)$ \\ \hline
$3$ & $83$ & $x^3+x^2+11$ & True \\
$4$ & $83$ & $x^4+x^3+24$ & True \\
\hline
\end{tabular}

\vspace{1em}

Primitive polynomials over $\FF_{89}$

\begin{tabular}{|c|c|c|c|c|} \hline
$n$ & $p$ & $P(x)$ & $B(n,p,P)$ \\ \hline
$3$ & $89$ & $x^3+x^2+6$ & True  \\
$4$ & $89$ & $x^4+x^3+14$ & True \\
\hline
\end{tabular}

\vspace{1em}

Primitive polynomials over $\FF_{97}$

\begin{tabular}{|c|c|c|c|c|} \hline
$n$ & $p$ & $P(x)$ & $B(n,p,P)$ \\ \hline
$3$ & $97$ & $x^3+x^2+5$ & True  \\
$4$ & $97$ & $x^4+x^3+15$ & True \\
\hline
\end{tabular}
\end{center}

\vspace{1em}

As we can see, our theorem can be applied in a lot of cases. But we could ask for an estimate of chances of success, that is to say: when $p \geq 5$ is a prime, and $n \geq 3$, what is the probability for $n$ to be a multiple of $p-1$? We hope that this probability is as small as possible.

If $f(X) \in \mathbb{F}_p[X]$ is a primtive polynomial of degree $n \geq 3$, of constant coefficient $a_0$, what is the probability that $n$ is not divisible by $p-1$ if $p \geq 5$ is prime?

For $N \geq 7$, let
\[
\mathcal{E}_N = \{ (p,n) \,\vert\, 5\leq p \leq N, p \in \mathcal{P}, 3\leq n \leq N \},
\]
\[
\mathcal{D}_N = \{ (p,n) \,\vert\, 5\leq p \leq N, p \in \mathcal{P}, 3\leq n \leq N, n \in (p-1)\mathbb{Z} \},
\]

We have: 
\[\mathcal{D}_N = \{ (p,k(p-1)) \,\vert\, 5\leq p \leq N, p \in \mathcal{P}, \frac{3}{p-1}\leq k \leq \frac{N}{p-1} \},
\]
\[
\# \mathcal{E}_N = (\pi(N)-2)(N-2). 
\]

The number of integers $k$ such that $\frac{3}{p-1}\leq k \leq \frac{N}{p-1}$ equals $\lfloor \frac{N}{p-1} \rfloor$, hence:
\[
\# \mathcal{D}_N = (\pi(N)-2) \lfloor \frac{N}{p-1} \rfloor. 
\]

Therefore:
\[
\frac{\mathcal{D}_N}{\mathcal{E}_N} = \frac{\lfloor \frac{N}{p-1} \rfloor}{N-2} \leq \frac{N}{(N-2)(p-1)} \leq \frac{N}{4(N-2)},
\]
because $p \geq 5$.

Thus:
\[
\lim_{N \to +\infty} \frac{\mathcal{D}_N}{\mathcal{E}_N}  \leq \frac{1}{4} \text{ and } \lim_{N \to +\infty} \frac{\mathcal{\overline{D}}_N}{\mathcal{E}_N} \geq\frac{3}{4}.
\]

If $n \geq 3$ and $p \geq 5$, the Brunnian conjecture is true in more than $75~\%$ of the cases.



\bibliography{amsbcjyd}{}
\bibliographystyle{alpha}

\end{document}